\documentclass[11pt]{article}

\usepackage{amscd,amsmath, amssymb, fancyhdr, graphicx, url}

\numberwithin{equation}{section}

\newcommand{\version}{version 1.0,\ \ July 30, 2019}

\makeatletter
\def\x@arrow{\DOTSB\Relbar}
\def\xlongrightarrowfill@{\arrowfill@\relbar\relbar\longrightarrow}
\newcommand{\xlongrightarrow}[2][]{%
        \ext@arrow 0099\xlongrightarrowfill@{#1}{#2}}
\makeatother

\def\eqref#1{(\ref{#1})}

\newcommand{\arrow}{{\:\longrightarrow\:}}
\newcommand{\Z}{{\Bbb Z}}
\def\C{{\Bbb C}}
\def\P{{\Bbb P}}

\newcommand{\R}{{\Bbb R}}
\newcommand{\Q}{{\Bbb Q}}

\def\1{\sqrt{-1}\:}

\newcommand{\cntrct}                
{\hspace{2pt}\raisebox{1pt}{\text{$\lrcorner$}}\hspace{2pt}}

\newcommand{\calo}{{\cal O}}

\renewcommand{\tilde}{\widetilde}
\renewcommand{\bar}{\overline}
\renewcommand{\phi}{\varphi}
\renewcommand{\epsilon}{\varepsilon}
\renewcommand{\geq}{\geqslant}
\renewcommand{\leq}{\leqslant}

\renewcommand{\min}{{\operatorname{\sf min}}}
\newcommand{\Teich}{\operatorname{\sf Teich}}
\newcommand{\Hilb}{\operatorname{Hilb}}
\newcommand{\Comp}{\operatorname{\sf Comp}}

\newcommand{\Per}{\operatorname{\sf Per}}
\newcommand{\Perspace}{\operatorname{{\Bbb P}\sf er}}

\newcommand{\Kah}{\operatorname{Kah}}

\newcommand{\Gr}{\operatorname{Gr}}

\newcommand{\Sym}{\operatorname{Sym}}
\newcommand{\Pic}{\operatorname{Pic}}
\newcommand{\Pos}{\operatorname{Pos}}

\newcommand{\Mon}{\operatorname{\sf Mon}}

\newcommand{\Diff}{\operatorname{\sf Diff}}

\renewcommand{\Re}{\operatorname{Re}}
\renewcommand{\Im}{\operatorname{Im}}

\newcommand{\proof}{\noindent{\bf Proof:\ }}

\newcounter{Mycounter}[section]
\newcounter{lemma}[section]
\renewcommand{\thelemma}{{Lemma \thesection.\arabic{lemma}}}
\newcommand{\lemma}{%
    \setcounter{lemma}{\value{Mycounter}}
    \refstepcounter{lemma}
    \stepcounter{Mycounter}
    {\noindent \bf \thelemma:\ }}

\newcounter{claim}[section]
\newcounter{sublemma}[section]
\newcounter{corollary}[section]
\renewcommand{\thecorollary}{{Corollary \thesection.\arabic{corollary}}}
\newcommand{\corollary}{%
    \setcounter{corollary}{\value{Mycounter}}
    \refstepcounter{corollary}
    \stepcounter{Mycounter}
    {\noindent \bf \thecorollary:\ }}

\newcounter{theorem}[section]
\renewcommand{\thetheorem}{{Theorem \thesection.\arabic{theorem}}}
\newcommand{\theorem}{%
    \setcounter{theorem}{\value{Mycounter}}
    \refstepcounter{theorem}
    \stepcounter{Mycounter}
    {\noindent \bf \thetheorem:\ }}

\newcounter{conjecture}[section]
\newcounter{proposition}[section]
\renewcommand{\theproposition}
      {{Proposition \thesection.\arabic{proposition}}}
\newcommand{\proposition}{%
    \setcounter{proposition}{\value{Mycounter}}
    \refstepcounter{proposition}
    \stepcounter{Mycounter}
    {\noindent \bf \theproposition:\ }}

\newcounter{definition}[section]
\renewcommand{\thedefinition}
      {{Definition~\thesection.\arabic{definition}}}
\newcommand{\definition}{%
    \setcounter{definition}{\value{Mycounter}}
    \refstepcounter{definition}
    \stepcounter{Mycounter}
    {\noindent \bf \thedefinition:\ }}

\newcounter{example}[section]
\renewcommand{\theexample}{{Example \thesection.\arabic{example}}}
\newcommand{\example}{%
    \setcounter{example}{\value{Mycounter}}
    \refstepcounter{example}
    \stepcounter{Mycounter}
    {\noindent \bf \theexample:\ }}

\newcounter{remark}[section]
\renewcommand{\theremark}{{Remark \thesection.\arabic{remark}}}
\newcommand{\remark}{%
    \setcounter{remark}{\value{Mycounter}}
    \refstepcounter{remark}
    \stepcounter{Mycounter}
    {\noindent \bf \theremark:\ }}

\newcounter{problem}[section]
\newcounter{question}[section]
\makeatletter

\@addtoreset{equation}{section} \@addtoreset{footnote}{section}
\makeatother

\def\blacksquare{\hbox{\vrule width 5pt height 5pt depth 0pt}}
\def\endproof{\blacksquare}

\begin{document}

\begin{center}
{\Large\bf
MBM classes and contraction loci on low-dimensional hyperk\"ahler manifolds of K3${}^{[n]}$ type\\[3mm]
}

Ekaterina Amerik\footnote{Partially supported 
by  the  Russian Academic Excellence Project '5-100'.}, 
Misha Verbitsky\footnote{Partially supported 
by  the  Russian Academic Excellence Project '5-100', 
FAPERJ E-26/202.912/2018 and CNPq - Process 313608/2017-2.


}

\end{center}

{\small \hspace{0.15\linewidth}
\begin{minipage}[t]{0.7\linewidth}
{\bf Abstract} \\
An MBM locus on a hyperk\"ahler manifold is the 
union of all deformations of a minimal rational
curve with negative self-intersection.
MBM loci can be equivalently defined as 
centers of bimeromorphic contractions. 
In \cite{_AV:contr_}, it was shown that 
the MBM loci on deformation equivalent
hyperk\"ahler manifolds are diffeomorphic.
We determine the MBM loci on a hyperk\"ahler
manifold of K3${}^{[n]}$-type for small $n$ using
a deformation to a Hilbert scheme of a
non-algebraic K3 surface. 
\end{minipage}
}

\tableofcontents


\section{Introduction}


In the last twenty years or so, a considerable effort has been made to understand the Mori cone, or, dually, 
the ample (or K\"ahler in the non-projective case) cone of irreducible holomorphic symplectic (IHS) manifolds. For K3 surfaces it has been known for a long time 
that the extremal rays correspond to $(-2)$-curves (that
is, irreducible curves of square $-2$, automatically
smooth and rational), so that the ample (or K\"ahler)
classes are characterized, among positive classes, as
those which are positive on all $(-2)$-curves. In higher
dimension, though, until recently the only general result,
due to Huybrechts and Boucksom, was that one could still
test the K\"ahler property on rational curves.
However, no boundedness statement for such curves was
available. A boundedness statement (equivalent to the
Morrison-Kawamata cone conjecture) was obtained in 
\cite{_AV:Mor_Kaw_}.

In \cite{HT1} Hassett and Tschinkel conjectured a
description of the Mori cone for projective IHS manifolds
which are deformation equivalent to the Hilbert square of
a K3 surface. Such manifolds are also called holomorphic
symplectic fourfolds of K3${}^{[n]}$ type.
The conjecture by Hassett and Tschinkel (Conjecture 3.1 of
\cite{HT1})  stated that the extremal rays are generated
by the  classes with some particular numerical
properties. For the Hilbert square of a K3 surface
containing a smooth rational curve $C$, these numerics are either
those of the rational curves contracted 
by the Hilbert-Chow morphism, or of the rational 
curves arising from $C$, or of the line in the projective
plane which is the symmetric square of such a $C$. Hassett
and Tschinkel have subsequently proved this conjecture in
\cite{HT2} (in the projective case). 

For high-dimensional holomorphic symplectic manifolds of K3${}^{[n]}$ type 
the analogous initial conjecture by Hassett and Tschinkel 
(Conjecture 1.2 of \cite{HT3} turned out to be false. The description of the Mori cone was provided by Bayer and 
Macri \cite{BM} using advanced technology around Bridgeland's stability conditions.

Another conjecture by Hassett and Tschinkel (Conjecture 3.6 of \cite{HT1}, and still a conjecture 24\ in \cite{HT2})  was concerned with the geometry of the loci covered by the minimal rational 
curves
according to the numerical properties of the corresponding
extremal rays (by Kawamata base-point-free theorem, one
can interpret those loci
as exceptional sets of birational contractions). As it is
natural, one expects that such a locus is either 
a lagrangian plane or a $\P^1$-bundle over a K3 surface, the former covered by curves with Beauville-Bogomolov square 
$-5/2$, the latter either $-1/2$ (the type of the
exceptional divisor of the Hilbert-Chow map) or $-2$ (the
type of the divisor of subschemes with some support on 
a $-2$-curve $C$ on the K3-surface). One can also expect a similar
description in higher dimension. Hassett and Tschinkel
have given lists of possible types of extremal rays and
corresponding contraction loci in \cite{HT3} and then
completed a few missing types following Bayer and Macri in
\cite{HT4}.

Bayer and Macri's technology allows, in principle, to describe contraction loci as fibrations in projective spaces
 over 
K3-related manifolds (up to birational equivalence). Those manifolds arise as moduli spaces of sheaves on a K3 surface in question, and the 
projective spaces are projectivisations of Ext-groups between the sheaves. But this approach is not very explicit, 
whereas at list in low dimension it seems possible and useful to describe the extremal rays and the subvarieties covered by the corresponding curves ``by hand''.

The purpose of this note is to provide a particularly simple description of possible extremal rays and contraction loci for holomorphic symplectic manifolds of K3${}^{[n]}$ type in low dimension, particularly four and six
(i.e. for manifolds which are deformation equivalent to the Hilbert square or the Hilbert cube of a K3 surface). Rather than studying extremal rays, we work with a slightly different notion of MBM classes. This notion was introduced in \cite{_AV:MBM_} to put things in the deformation-invariant context, and indeed the crucial point in our argument is the deformation to the Hilbert scheme of a non-projective K3 surface.

 The MBM classes are the classes which correspond to extremal rays of the Mori cone on some deformation of $X$, or, equivalently, the classes of extremal rays up to monodromy and birational equivalence (the name derives from ``monodromy birationally minimal''). From a uniruled 
subvariety $Z$ of an IHS $X$ we obtain an MBM class as follows: let $C$ be a rational curve 
which is minimal among the rational curves in $Z$ (in the sense that $C$ does not bend-and-break). If its class is Beauville-Bogomolov negative, then it is MBM. 

The parameter space $\Teich_{\alpha}^{min}$ for complex structures where an MBM class $\alpha$ is represented by a minimal rational 
curve coincides with the whole space where $\alpha$ stays
of type $(1,1)$ up to inseparability issues. Such curves
can always be contracted (\cite[Corollary 5.3]{_AV:contr_}).

In \cite{_AV:contr_} we proved some deformation-invariance statements about the corresponding contraction 
loci (i.e. exceptional sets of morphisms contracting the class $\alpha$). We have, in particular, shown that the diffeomorphism type of such a contraction locus is constant over $\Teich_{\alpha}^{min}$. Moreover, so is the biholomorphism type of the fibers of the rational quotient fibration of the exceptional locus, with a possible exception of the points of $\Teich_{\alpha}^{min}$ corresponding to the complex structures with maximal Picard number.

The idea now is as follows. Let $\alpha$ be a negative integral $(1,1)$-class on a hyperk\"ahler fourfold $X$ of K3${}^{[n]}$ type. Using Markman's description of the monodromy group, we prove that $X$ can be deformed to a non-projective manifold $X'\cong \Hilb^2(S')$ where $S'$ is a K3 surface with cyclic Picard group generated by a 
negative class, and in such a way that  $\alpha$ stays of type $(1,1)$. Since rational curves on $X'$ are easy to describe, we obtain the classification of subvarieties which can be covered by rational curves in the class $\alpha$, by deformation-invariance.

If $X$ is an IHS of K3${}^{[n]}$ type of dimension six, it
turns out that $X$ still deforms to the Hilbert cube of a
non-projective K3 surface $S$, but its Picard group can
now be generated by an element with negative or zero
Beauville-Bogomolov square. The geometry in this last case
is also well-understood. Namely, $S$ is an elliptic $K3$
surface with singular fibers which are nodal rational
curves, and there is no other curve at all. All rational
curves (up to deformation)
on the symmetric cube of $S$ are obtained from linear
systems  $g_k^1, k\leq 3$ on the fibers,\footnote{As
  defined in \cite[\S IV.5]{_Hartshorne_}, $g_k^d$
denotes a linear system of dimension $d$ and degree $k$
on a curve. In particular, $g_k^1$ correspond to
$k$-sheeted maps to $\P^1$.}
so again easy to
describe. We obtain a table of five types of MBM classes,
their squares and the types of the corresponding
contraction loci, which agrees with the data by Hassett
and Tschinkel (one may notice here that the MBM class
initially overlooked in \cite{HT3} and exhibited in
\cite{HT4} while working out Bayer and Macri's formulae on
this example, comes in our picture from the Hilbert cube
of a surface with zero intersection form on the Picard
group).

The same calculation works in fact for eightfolds (again bringing into the picture a class originally overlooked 
in \cite{HT3}) and tenfolds, but not for higher-dimensional IHS of K3${}^{[n]}$ type. It would be interesting to see 
whether the present elementary method could be adapted to give a reasonable classification in all dimensions.


\section{Hyperk\"ahler manifolds}


In this section we recall some definitions concerning IHS manifolds and rational curves.
For more details see
\cite{_Beauville_} and \cite{_AV:MBM_}.

\subsection{Holomorphically symplectic manifolds}

\definition {\bf A holomorphically symplectic manifold}
is a complex manifold equipped with a nowhere degenerate holomorphic
two-form.

\hfill

\definition For the rest of this paper,
{\bf a hyperk\"ahler manifold}
is a compact, K\"ahler, holomorphically symplectic manifold.

\hfill

\definition A hyperk\"ahler manifold $M$ is called  {\bf
  simple}, or {\bf IHS},   or {\bf of maximal holonomy}
if $\pi_1(M)=0$, $H^{2,0}(M)=\C$.

\hfill

The terminology originates from {\bf Bogomolov's decomposition theorem} stating that 
any hyperk\"ahler manifold admits a finite covering which 
is a product of a torus and several  IHS,
\cite{_Bogomolov:decompo_}. This decomposition also
corresponds the de Rham-Berger local decomposition 
of a Riemannian manifold into a product of manifolds
with irreducible Riemannian holonomy. Then, the IHS manifolds correspond
to the hyperk\"ahler manifolds which have irreducible holonomy.
Further on, all hyperk\"ahler manifolds
are assumed to be of maximal holonomy.

\hfill

\example The basic series of examples is that of punctual Hilbert schemes of K3 surfaces. As observed by 
Beauville in \cite{_Beauville_}, for any K3 surface $S$, the Hilbert scheme $S^{[n]}$ of ideal sheaves of colength $n$ 
on $S$ is an IHS. An IHS which is a deformation of $S^{[n]}$ for a K3 surface $S$ is said to be {\bf of K3${}^{[n]}$ type}.   There is one more series (generalized Kummer manifolds) and two sporadic families (O'Grady's examples) 
known, and we don't know whether there is anything else.

\subsection{The Bogomolov-Beauville-Fujiki form}

\theorem (Fujiki, \cite{_Fujiki:HK_}) 
Let $\eta\in H^2(M)$, and $\dim M=2n$, where $M$ is
hyperk\"ahler. Then $\int_M \eta^{2n}=cq(\eta,\eta)^n$,
for some primitive integral quadratic form $q$ on $H^2(M,\Z)$,
and $c>0$ a rational number.

\hfill

\definition This form is called
{\bf Bogomolov-Beauville-Fujiki form}. It is defined
by the Fujiki's relation uniquely, up to a sign. The sign is found
from an explicit formula due to Beauville and Bogomolov, using integration against exterior powers of the holomorphic symplectic form $\Omega$ and its conjugate, see e.g. \cite{_Beauville_}.

\hfill

\remark  The form $q$ has signature $(3,b_2-3)$.
It is negative definite on primitive forms, and positive
definite on $\langle \Omega, \bar \Omega, \omega\rangle$,
 where $\omega$ is a K\"ahler form. 

\remark The BBF form is in general not unimodular. We shall use it to embed
$H_2(M,\Z)$ into $H^2(M,\Q)$ as an overlattice of $H^2(M,\Z)$. This allows us to view the classes of curves on $M$
as rational $(1,1)$-classes and compute the BBF square of a curve which is a rational number (cf. \cite{HT1}).

\hfill

\example\label{bbf-K3}If $M$ is a K3 surface, the BBF form is just the intersection form. If $M=S^{[n]}$ for a K3 surface $S$,
we have an orthogonal direct sum decomposition of lattices $H^2(M, \Z)=H^2(S,\Z)\oplus \Z e$, where $H^2(S,\Z)$ is equipped with the intersection form and  $e^2=-2(n-1)$. In fact $e$ is one half of the class of the diagonal, that is, the divisor of non-reduced subschemes of length $n$ in $S$ (this divisor is also the exceptional set of the Hilbert-Chow morphism for $S^{[n]}$.
This decomposition holds for IHS of K3${}^{[n]}$ type (but the class $e$ may cease to be algebraic).

\subsection{Extremal curves}

Let $Z\subset M$ be an irreducible uniruled subvariety of a compact K\"ahler manifold $(M, \omega)$. 
We call a rational curve $C\subset Z$ minimal if it is of minimal degree with respect to the K\"ahler form, 
i.e. $\int_{C} \omega$ is minimal, among the curves passing through a general point of $Z$.


\hfill

\remark Thanks to Mori's bend-and-break lemma, there is the following necessary condition
for minimality: for general $x,y\in Z$, the space of deformations
$Z_{x,y}$ of $C$ passing through $x,y$ is 0-dimensional. Indeed, if  $\dim Z_{x,y}>0$, we can use Mori lemma to deform $C$  to a reducible curve which contains $x,y$. The area $\int_{C'} \omega$ of each
component of this union is strictly smaller than that of $C$.

\hfill

If $X$ is IHS, minimal rational curves on $X$ have the following remarkable properties (see for instance \cite{_AV:MBM_}).

\hfill

\proposition\label{min-def} Let $X$ be IHS of dimension $2n$, then a minimal rational curve $C$ (with respect to some maximal uniruled subvariety $Z$) deforms together with $X$ whenever its cohomology class remains of type $(1,1)$.
The space of deformations of $C$ within $X$ is of dimension $2n-2$. The subvariety $Z$ is coisotropic and can be of any dimension between $n$ and $2n-1$. Its codimension in $X$ is equal to the relative dimension of its rational quotient fibration. 

\hfill



\definition
A rational curve is called {\bf MBM} if it is minimal and
its cohomology class has negative self-intersection.

\hfill

One reason why we are especially interested in negative
self-intersection rational curves is that they determine
the K\"ahler cone $\Kah(M)$ (or the ample cone in the
projective case) inside the 
{\bf positive cone}. The latter is defined as the connected component of the set 
\begin{equation}\label{_pos_cone_Equation_}
\{\eta \in H^{1,1}(M,\R)\ \ |\ \ (\eta,\eta)>0\}
\end{equation}
 which contains the K\"ahler cone $\Kah(M)$. Note that the BBF form on $H^{1,1}(M,\R)$ has signature $(1, b_2-3)$
and therefore the set \eqref{_pos_cone_Equation_} 
has two connected components.


\hfill

\theorem (Huybrechts, Boucksom, see \cite{_Boucksom-cone_}) 
The K\"ahler cone of $M$ is the set of all $\eta\in \Pos(M)$
such that $(\eta, C)>0$ for all rational curves $C$.


\hfill

Clearly we may consider only negative rational curves, as by the Hodge index theorem the positive classes never have zero intersection with curves of positive square. If moreover the class of such a curve cannot be written as a sum
of two other such classes, it is called an {\bf extremal} rational curve. The K\"ahler cone is locally polyhedral (with
some round pieces in the boundary), and its faces are
orthogonal complements to the classes of extremal curves. Note that an MBM curve as in \ref{min-def} is extremal on a 
generic deformation, as only its class, up to proporionality, survives as a $(1,1)$-class there.

\subsection{MBM classes}

A minimal rational curve may split, and thus cease to be
minimal, under a deformation. To put the things into a
deformation-invariant context, we have introduced the
notion of {\bf MBM classes} in \cite{_AV:MBM_}. These are,
up to a rational multiple, monodromy transforms (see next
subsection for the definition of the monodromy) of the classes of
extremal rational curves on IHS birational models of $M$. 
Note that if $M$ and $M'$ are two birational IHS, or, more
generally, two birational manifolds with trivial canonical
bundle, then $H^2(M)$ is naturally identified with
$H^2(M')$. The following is a reformulation of 
\cite[Theorem 1.17]{_AV:MBM_}.


\hfill

\theorem 
Let $\eta$ be the cohomology class of an extremal curve on a 
hyperk\"ahler manifold, and
$(M_1,\eta)$ be obtained as a deformation of $(M,\eta)$
in a family such that $\eta$ remains of type $(1,1)$
for all fibers of this family. Then $M_1$ is birational to
a hyperk\"ahler manifold 
$M_1'$ such that $\gamma(\eta)$ is a class of an extremal
curve on $M_1'$ for some $\gamma$ in the monodromy group.



\hfill

\remark Birational hyperk\"ahler manifolds are deformation equivalent by a result of Huybrechts. In fact a more precise version of the theorem above
(see \cite{_AV:contr_}) states that an extremal rational curve remains extremal on an $M_1'$ unseparable from $M_1$ in the 
Teichm\"uller space (see next section for preliminaries on Teichm\"uller space and more details).

\hfill

MBM classes can be defined in several equivalent ways. The simplest  one to formulate is as follows.

\hfill

\definition
A class $\eta\in H^2(M,\Z)$ which, up to a rational multiple, can be represented by an extremal curve
on a deformation on $M$
is called {\bf an MBM class}.

\hfill

\remark Another equivalent definition (\cite{_AV:MBM_}):
An MBM class $\eta\in H^2(M)$ is a Beauville-Bogomolov negative class such that $\lambda\eta$ can be represented by a rational curve in some $(M,I)$ when $(M,I)$ is a non-algebraic deformation of $M$ with $\Pic(M, I)_\Q =\langle \eta\rangle$.
In that case it can be represented by a rational curve in all such $(M,I)$.

\hfill

The set of MBM classes is thus the same for a given deformation type. If a complex structure $I$ is fixed, some of them are of type $(1,1)$, and among those, some correspond to extremal curves (that is, determine the K\"ahler cone), others correspond to extremal curves on birational models and their monodromy transforms. The following result further explains the connection between MBM classes and the K\"ahler cone.  

\hfill








\theorem (\cite{_AV:MBM_}
Let $(M,I)$ be a hyperk\"ahler manifold and $S$ the set of all
MBM classes of type (1,1) on $(M,I)$. Let $S^\bot$ the union of all
orthogonal complement to all $s\in S$. Then $\Kah(M,I)$ is 
a connected component of $\Pos(M)\backslash S^\bot$. 

\hfill

Following E. Markman, we call these connected components
{\bf the K\"ahler-Weil chambers}.


\section{Teichm\"uller spaces, Torelli theorem and monodromy}



\subsection{Teichm\"uller space}

This subsection and the next one follow \cite{_V:Torelli_}.

\hfill

\definition
Let $M$ be a hyperk\"ahler manifold, and 
$\Diff_0(M)$ a connected component of its diffeomorphism group
(the group of isotopies). Denote by $\Comp$
the space of complex structures of K\"ahler type on $M$, and let
$\Teich:=\Comp/\Diff_0(M)$. We call 
it {\bf the Teichm\"uller space.}

\hfill

\remark The space $\Teich$ is a finite-dimensional
complex manifold (\cite{_Catanese:moduli_}), not necessarily Hausdorff.

\hfill

\definition
Let $\Diff(M)$ be the group of 
diffeomorphisms of $M$. We call $\Gamma:=\Diff(M)/\Diff_0(M)$ {\bf the
mapping class group}. 

\hfill

The quotient $\Teich/\Gamma$ 
is identified with the set of equivalence classes of complex structures. However the notion of the moduli space usually does not make sense as the action of $\Gamma$ on $\Teich$ is ``too bad''.

By a result of Huybrechts, the quotient 
$\Comp/\Diff_0(M)$ has only finitely many connected components. To simplify notations we denote by $\Teich$ the one which contains our given complex structure, that is, a parameter space for the deformations of $M$. The group $\Gamma$ shall mean the subgroup of the mapping class group which preserves $\Teich$. 

\hfill

\definition 
The {\bf monodromy group} $\Mon$ is the image of the representation of $\Gamma$
on $H^2(M, \Z)$. 
 
\hfill

\definition Let 
$\Per:\; \Teich \arrow {\Bbb P}H^2(M, \C)$
map a complex structure $J$ to the line $H^{2,0}(M,J)\in {\Bbb P}H^2(M, \C)$.
The map $\Per:\; \Teich \arrow {\Bbb P}H^2(M, \C)$ is 
called {\bf the period map}. The map $\Per$ maps $\Teich$ to an open subset of a 
quadric
\[
\Perspace:= \{l\in {\Bbb P}H^2(M, \C)\ \ | \ \  q(l,l)=0, q(l, \bar l) >0\}.
\]
It is called {\bf the period space} of $M$.

\hfill

\remark 
The space $\Perspace$ can be naturally
identified with 
\[ SO(b_2-3,3)/SO(2) \times
SO(b_2-3,1)=\Gr_{++}(H^2(M,\R))
\] 
(the Grassmannian of Beauville-Bogomolov positive, oriented planes in $H^2(M,\R)$).

Indeed, the group $SO(H^2(M,\R),q)=SO(b_2-3,3)$ acts transitively on
$\Perspace$, and $SO(2) \times SO(b_2-3,1)$ is the
stabilizer of a point. More explicitly, we can associate
to each complex line a real plane generated by its real
and imaginary part.

\subsection{Global Torelli theorem}

\definition
Let $M$ be a topological space. We say that $x, y \in M$
are {\bf non-separable} (denoted by $x\sim y$)
if for any open sets $V\ni x, U\ni y$, $U \cap V\neq \emptyset$.

\hfill






\definition
The space $\Teich_b:= \Teich/\sim$ is called {\bf the
birational Teichm\"uller space} of $M$. Since $\Teich_b$
is obtained by gluing together all non-separable points,
it is also called {\bf the Hausdorff reduction} of $\Teich$.

\hfill

\theorem (Global Torelli theorem for hyperk\"ahler manifolds)\\
The period map 
$\Teich_b\stackrel \Per \arrow \Perspace$ is a diffeomorphism,
for each connected component of $\Teich_b$.

\hfill

Huybrechts in \cite{_Huybrechts:basic_} observed that unseparable points of $\Teich$ correspond to birational IHS,
up to the action of the mapping class group. The following precise description of unseparable points is due to Markman.

\hfill

\theorem \label{_non_Hausdorff_points_Markman_Theorem_}
The points in $\Per^{-1}(\Per I)$ correspond to the K\"ahler-Weil chambers (that is, to the connected components
of $\Pos(M)\backslash \cup z^{\bot}$ where $z$ runs through the MBM classes of type $(1,1)$).
Each chamber is the K\"ahler cone of the corresponding point in $\Per^{-1}(\Per I)$.

\hfill

In particular, as a general hyperk\"ahler manifold has no curves (the ones which have curves belong to a countable union of divisors in 
$\Teich$), $\Teich$ is generically Hausdorff.


\subsection{The space $\Teich_{\alpha}^{\min}$}

Let now $\alpha$ be an MBM class, we are interested in the complex structures $I$ from $\Teich$ such that $\alpha$
is of type $(1,1)$ in $I$. Denote the subspace of such
structures in $\Teich$ by $\Teich_{\alpha}$.
Clearly, $\Teich_{\alpha}=\Per^{-1}(\alpha^{\bot})$. 
From \ref{_non_Hausdorff_points_Markman_Theorem_} it is clear that $\Teich_{\alpha}$
is not Hausdorff even generically (indeed each complex structure $I$ in $\Teich_{\alpha}$ has at least two 
K\"ahler-Weil chambers in which $\alpha^{\bot}\subset H^{1,1}(I)$ divides the positive cone). On the other hand, 
$\Teich_{\alpha}$ is naturally divided in two halves, according to whether ${\alpha}$ is negative or positive on the
K\"ahler classes. In the ``positive'' half $\Teich_{\alpha}^+$ take the part formed by the complex structures whose 
K\"ahler cones are adjacent to ${\alpha}^\bot$.

\hfill

\definition The space thus defined is called $\Teich_{\alpha}^{\min}$. In other words $\Teich_{\alpha}^{\min}$ parameterizes the complex 
structures such that
${\alpha}$ is the class of a minimal rational curve, up to a positive multiple.

\hfill

\remark It is clear from this description that $\Teich_{\alpha}^{\min}$ is generically Hausdorff and that $\Per$
is an isomorphism from $\Teich_{\alpha, b}^{\min}$ to $\alpha^{\bot}\subset\Perspace$.

\subsection{Monodromy for manifolds of K3${}^{[n]}$ type}

The monodromy group is very important for understanding the geometry of IHS manifolds, but in the general (``classification-free'') case its structure remains somewhat mysterious. One knows, thanks to Markman \cite{Markman-pex}, that the
reflections in the classes of uniruled divisors are in $\Mon$, and this already has nontrivial consequences. But it is not clear how to decide whether a given MBM class
is ``divisorial'' or not (that is, whether the corresponding rational curves cover a divisor or a subvariety of higher codimension). 

A lot of work has been done on the known examples of
IHS. For manifolds of K3${}^{[n]}$ type, Markman proved the
following.

\hfill

\theorem (\cite[Theorem 1.2]{_Markman:constra_})
The monodromy group $\Mon$ is the subgroup of 
$O(H^2(M,\Z))$ generated by the reflections in the classes of
square $-2$ and and the antireflections in those of square
$2$, that is, by the transformations
$\rho_u(x)=\frac{-2x}{q(u)}+q(x,u)u$, where the
Beauville-Bogomolov square
$q(u)=\pm 2$, and the same notation $q$ is used for the
 bilinear form associated to $q$. Up to $\pm$ sign,
it also coincides with the {\bf stable orthogonal group}
$\tilde{O}(H^2(M, \Z))$ of isometries acting trivially on
the discriminant group.

\hfill

Here the {\bf discriminant group} of a lattice $L$ is the quotient $L^*/L$ where $L^*$ is the dual lattice and
$L$ is naturally embedded into $L^*$. If $L$ is $H^2(M, \Z)$ for $M$ a $2n$-dimensional IHS of K3${}^{[n]}$ type,
the discriminant group is cyclic of order $2n-2$ (see \ref{bbf-K3}). We shall often view it as a subgroup of $\Q/\Z$, like in \ref{orb-classif} just below.

\hfill

\remark When $n-1$ is a power of a prime number, 
one can show that $\Mon$ coincides with the orientation-preserving part of the 
orthogonal group itself (\cite{_Markman:constra_}, lemma 4.2). We shall
not need this as the group $\tilde{O}$ has 
also been studied. Notably we shall use the 
following.

\hfill

\theorem (see e.g. \cite{GHS}, Lemma 3.5)\label{orb-classif} Let $(L, q)$ be a lattice containing two orthogonal isotropic planes. The orbit of $l\in L$
under $\tilde{O}(L)$ is determined by $q(l)$ and the image of $l/d(l)$ in the discriminant
 group. Here we view the dual lattice as an overlattice
of $L$ in $L\otimes \Q$, and $d(l)$ is the {\it divisibility} of $l$, that is, the pairing
 of $l$ with $L$ gives the ideal $d(l)\Z$.


\section{Birational contractions and invariance by deformation}
\label{_BPF_Section_}

\subsection{Birational contractions}

In this section we recall the Kawamata base point free theorem which implies that extremal rational curves
on projective IHS manifolds are contractible, as well as a more recent observation by Bakker and Lehn concerning 
the nonprojective case.



\hfill

\definition
A line bundle $L$ on a compact complex k-fold $M$ is {\bf nef} if $c_1(L)$ lies in the closure of the
K\"ahler cone, {\bf big} if $H^0(M, L^{\otimes N})= O(k^N)$ for $N$ sufficiently big and sufficiently divisible, 
and  {\bf semiample} if $L^{\otimes N}$ is generated by its global sections for such $N$.

\hfill

\theorem (Kawamata base point free theorem, \cite[Theorem 4.3]{_Kawamata:Pluricanonical_})\\
Let $L$ be a nef line bundle on a projective $M$ such that $L^{\otimes N}\otimes {\cal O}(-K_M)$ is nef and big for some $N$.
Then $L$ is semiample.


\hfill

For manifolds with trivial canonical bundle this just means that big and nef bundles are semiample.




It follows that the sections of a suitable power of a big
and nef bundle $L$ on such manifolds define a
birational contraction (that is, a birational morphism with connected fibers) $\phi:M \arrow Y$.

Suppose now that $M$ is a projective IHS manifold,
$\alpha$ is an MBM class and the given complex structure
$I$ on $M$ is such that
$I\in \Teich_{\alpha}^{\min}$, so that $\alpha$ is in fact the class of an extremal curve: the hyperplane $\alpha^{\bot}$ in 
$H^{1,1}(M,\R)$ contains a
wall of the K\"ahler cone. Since $M$ is projective, the interior of this wall contains an integral $(1,1)$-class which by Lefschetz theorem 
is $c_1(L)$ for a line bundle $L$. This line bundle is nef, and also big. Indeed its BBF square is positive, and hence so is $c_1(L)^{2n}$.
Therefore a power of $L$ defines a birational morphism contracting exactly the curves of class proportional to $\alpha$.

It turns out that extremal curves can be contracted also in the non-projective case. This is essentially due to Bakker and Lehn in \cite{BL}. More precisely, putting together some results in \cite{_AV:MBM_}, \cite{BL} and \cite{_AV:contr_} gives the following.

\hfill

\theorem (\cite[Proposition 4.10]{BL}, \cite[Corollary
  5.3]{_AV:contr_}) \\
Let $M$ be an IHS manifold and $\alpha$
the class of an extremal rational curve on $M$. Then there
exists a birational contraction $\phi:M \arrow Y$
contracting to points the curves whose class is
proportional to $\alpha$. Moreover, for every component $Z$ of the
exceptional set, the fibers of $\phi|_Z$ are rationally
connected, and $\phi|_Z$ is the rational quotient
fibration of $Z$.

\hfill

In particular the uniruled subvarieties of $M$ covered by deformations of an extremal rational curve are irreducible components of the
exceptional sets of birational contractions. We sometimes refer to these as {\bf MBM loci}.

\subsection{Invariance by deformation}

The following is the main result of
\cite{_AV:contr_}.

\hfill


\theorem \label{_homeomo_loci_Theorem_} Let $\alpha$ be an MBM class of type $(1,1)$ 
on a hyperk\"ahler manifold $M$ with $b_2(M)>5$. Assume that $M$ and its deformation $M'$ are both in $\Teich_{\alpha}^{\min}$, so that
$\alpha$ is a class of an extremal rational curve on both.
Then the corresponding contraction centers on $M$ and $M'$ are sent one to another by a diffeomorphism $\Psi:M\to M'$ preserving the fibers of the contraction. If moreover neither $M$ nor
$M_1$ has maximal Picard rank, $\Psi$ can be chosen identifying the contracted extremal curves and thus induces the biholomorphism of the fibers provided that they are normal.

\hfill

\remark
Clearly, $H^{1,1}(M)$ is obtained as orthogonal complement
to the 2-dimensional space $\langle \Re\Omega, \Im \Omega\rangle$,
where $\Omega$ is the cohomology class of the holomorphic
symplectic form. Then $\Pic(M) = H^{1,1}(M)\cap H^2(M,\Z)$
has maximal rank only if the plane $\langle \Re\Omega, \Im \Omega\rangle$
is rational.  There is at most countable number of such $M$.

\hfill

\remark Existence of $\Psi$ follows from ergodicity of mapping class
group action, global Torelli theorem, Bakker-Lehn
contractibility result, and Thom-Mathers stratification
of proper real analytic maps.





\section{MBM loci for low-dimensional manifolds of $K3^{[n]}$ type}


The diffeomorphism theorem allows one to give an
explicit description of birational contraction centers and MBM loci in terms
of the period spaces and lattices. The easiest illustration is the (well-known, cf. Hassett and Tschinkel's works, 
but seemingly with no direct proof ever written) description
of the MBM loci on an IHS fourfold of K3${}^{[2]}$  type.

\hfill

\theorem\label{_K3^2_MBM_Theorem_} On an IHS fourfold of K3${}^{[2]}$  type, there are only three types of 
extremal rational curves $l$ (notation: $\alpha=[l]$, $Z$ is the MBM locus of $l$):
\begin{description}
\item[(a)] $q(\alpha)=-5/2$, $2\alpha$ is integral, $Z\cong\P^2$
\item[(b)] $q(\alpha)=-2$, $\alpha$ is integral, $Z$ is a $\P^1$-bundle over a K3 surface;
\item[(c)] $q(\alpha)=-1/2$, $2\alpha$ is integral, $Z$ is a $\P^1$-bundle over a K3 surface.
\end{description}

\hfill

Taking into account \ref{_homeomo_loci_Theorem_}, this is a consequence of the following precision:

\hfill

\theorem\label{_K3^2_MBM_Theorem_defo_} Let $\alpha$ be the class of an extremal rational curve on an IHS fourfold $M$ of K3${}^{[2]}$ type
and $Z$ the corresponding MBM locus. Then there exists a
deformation $M'$ of $M$ within $\Teich_{\alpha}$,
isomorphic
to the Hilbert square of a non-projective K3 surface $S$ with cyclic Picard group generated by a negative class $z$, such that one of the following holds:
\begin{description}

\item[(a)] $S$ contains a smooth rational curve $C$ and
  $Z$ deforms  to its Hilbert square, isomorphic to $\P^2$; $\alpha$ is 
a line in this $\P^2$;

\item[(b)] $S$ contains a smooth rational curve $C$ and $Z$ deforms to $Z_C\subset M'$, the divisor of all length-two subschemes of 
$S$ with support intersecting $C$; $Z_C$ is a ruled divisor singular along $C^{[2]}\cong \P^2$ and $\alpha$ is the class of the ruling.

\item[(c)] $Z$ deforms to the exceptional divisor of the
  Hilbert-Chow map \\ $S^{[2]}\to \Sym^2(S)$,
and $\alpha$ is the class of the ruling.
\end{description}

\hfill

\remark In the second case $\alpha$ is clearly not a class of an extremal rational curve, that is, $M'$ is not in $\Teich_{\alpha}^{min}$. As we know, this means that some complex structure unseparable from $M$ is in $\Teich_{\alpha}^{min}$. In fact, to render $\alpha$ extremal, one performs a flop in $C^{[2]}$, which also resolves the singularities of $Z_C$ and makes it a $\P^1$-bundle.

\hfill

Once the existence of a deformation $M'\cong S^{[2]}$, $\Pic(S)=\langle z\rangle$, $q(z)<0$, is shown, the rest is easy.
Indeed one has the following description of rational curves on $M'$.

\hfill

\proposition\label{rc-on-hilb-negative}
Let $S$ be a K3 surface with the Picard group $\Pic(S)$ of rank 1 generated
by a vector $z\in H^{1,1}(S,\Z)$ with negative square, and
$Z$ an MBM locus on $M^{[2]}$. Then $Z$ is described by
(a), (b) or (c) of \ref{_K3^2_MBM_Theorem_defo_}.

\hfill

\proof
If $S$ has no curves, the symmetric square $\Sym ^2 S$
also has no curves, and the bimeromorphic contraction
$S^{[2]}\arrow \Sym ^2 S$ contracts the exceptional
divisor $E$ to the singular set of $\Sym ^2 S$.
This means that $E=Z$, case (b) of
\ref{_K3^2_MBM_Theorem_defo_}.

If $S$ has a curve, it is a smooth (-2)-curve $C$ which
is necessarily unique. Indeed, any curve with negative self-intersection is a smooth rational curve with square $-2$ by adjunction, 
and there cannot be two different curves
in the proportional classes, again because $(z,z)<0$.
Then the (-2)-curve $C$ can be blown down to an ordinary double point.
Let $S_1$ be the result of this blow-down. Then 
$\Sym^2(S_1)$ has no curves, and the natural map
$S^{[2]}\arrow \Sym^2(S_1)$ blows down all curves
in $S^{[2]}$. The corresponding exceptional set $V$
is a union of two divisors, the exceptional divisor $E$
and the set $W\subset S^{[2]}$ of all ideal sheaves of colength two
$I\subset \calo_S$ with support of $I$ intersecting $C$.
Clearly, these two divisors intersect in $P^2=C^{[2]}$,
and outside of $C^{[2]}$ each of these divisors
is covered by a unique family of rational curves, hence the classification.
\endproof

\hfill

As for Beauville-Bogomolov squares, we recall the computations by Hassett and Tschinkel and others for convenience. 
We have the decomposition into the orthogonal direct sum 
$H^2(M,\Z)=H^2(S, \Z)\oplus \Z e$ where $e$ is half of the class of the diagonal $E$, and $e^2=-2$. 

\begin{itemize}

\item If $\alpha$ is the class of the ruling of the exceptional divisor, $q(e,\alpha)=-1$ since the exceptional divisor restricts as 
${\cal O}(-2)$ to the ruling. Hence $\alpha=e/2$.

\item If $\alpha$ is the class of the ruling of $Z_C$, then $\alpha$ is the class of $C$.

\item If $\alpha$ is the class of the line in the
  Lagrangian plane $C^{[2]}$, then $\alpha = [C]-e/2$
  (\cite{HT3}, Example 4.11; for $S^{[n]}$ this would be
  the same, indeed $q(e,\alpha)=n-1$)

\end{itemize}

We postpone the existence of $M'$ until the next section and deal with the next example of an IHS sixfold of K3${}^{[3]}$  type. 
There are five types of extremal rational curves/MBM classes/contraction loci. Four of them appear in the tables of \cite{HT3}
and are described geometrically exactly in the same way as those on a fourfold: on such an $M$ which is the Hilbert cube of a K3 
surface $S$ containing a (-2)-curve $C$, we have four natural subvarieties covered by deformations of rational curves in BBF-negative classes, that is, the exceptional divisor of the Hilbert-Chow map (ruled over the big diagonal), the ruled divisor of subschemes with some support on $C$ (ruled over the symmetric square of $S$), the fourfold
of subschemes with at least length-two part supported on $C$ (birational to a $\P^2$-bundle over $S$), and $\P^3=C^{[3]}$.

After Bayer and Macri have given their description of the Mori cone in \cite{BM}, Hassett and Tschinkel observed in \cite{HT4} that 
one case was actually missing from their table in dimension 6. The cohomology class of the extremal curve in question is half-integer 
and has square $-1$. We would like to present the following simple way to find this MBM class.

\hfill

\theorem\label{_K3^3_MBM_Theorem_defo_} Let $\alpha$ be a class of an extremal rational curve on an IHS $M$ of K3${}^{[n]}$ type of dimension
at most ten
and let $Z$ be the corresponding MBM locus. Then there exists a deformation $M'$ of $M$ within $\Teich_{\alpha}$, isomorphic
to the Hilbert scheme of a non-projective K3 surface $S$ with cyclic Picard group (thus generated by a non-positive class $z$).

\hfill

Again, let us describe the geometry we get in dimension six before going to the proof.

{\bf First case:} If $\Pic S$ is generated by a negative class, the rational curves on $M'$ are described as in \ref{rc-on-hilb-negative}:
one has the four types listed above. 

{\bf Second case:} If $\Pic S$ is generated by a class of square zero, then one deduces from Riemann-Roch theorem that 
$S$ carries an elliptic pencil with some degenerate fibers which are rational curves with a singularity, and no other curves.
For Hilbert cube this means that the extremal rational curves arise from linear systems of degree $k\leq 3$ on the (eventually normalized) curves from
the elliptic pencil: one fixes $3-k$ points on $S$ and lets the remaining $k$-uplet run through a pencil in the system. Of course, not all curves arising in this way are extremal, and some of them are already on the list obtained in the negative case.

Note that the computation of the BBF square is particularly easy since only the intersection with $e$ counts. Indeed
the intersection form is zero on its orthogonal complement
(which is  naturally identified with $\Pic S$).

For instance, the curve obtained by fixing two points and letting the third one run through the singular fiber is of BBF square 
zero, so not extremal. Considering a rational curve $R$
coming from a $g^1_3$ on a smooth elliptic curve
, we see that $q(R)=-9/4$ and $R$ is quarter-integer ($e$
comes with coefficient $-3/4$), by the following
well-known lemma (see \cite{CK} for a version which also 
applies to singular curves).

\hfill

\lemma\label{fromcoverings} The class of a rational curve coming from a $g^1_k$ on a curve $C$ of genus $g$ is 
$C-\frac{g+k-1}{2(n-1)}e$.

\hfill

The idea is that the intersection with the diagonal 
is computed by the Hurwitz formula, as this is roughly speaking the ramification of the morphism to the projective line 
defined by the $g^1_3$.

Such an extremal class (quarter-integer with square $-9/4$) is already obtained by Hassett-Tschinkel's computation when $S$ contains
a $(-2)$-curve $C$. 
It appears as the class of a line in the projective plane $C^{[2]}$. By the description of monodromy orbits, the two classes
are conjugate by a monodromy transform, so we get no new geometry here.

On the contrary, fixing one point on $S$ and considering a $g^1_2$ on an elliptic curve, one gets a half-integer class of square 
$-1$, exactly the one omitted in the table of \cite{HT3}.

The following corollary results from the above computations.

\hfill

\corollary\label{_K3^3_MBM_Theorem_} There are 5 types of MBM classes on an IHS sixfold of K3${}^{[3]}$  type:

\begin{description}
\item[(a)] $q(\alpha)=-3$, $2\alpha$ is integral, $Z\cong\P^3$
\item[(b)] $q(\alpha)=-9/4$, $4\alpha$ is integral, $Z$ is birationally a $\P^2$-bundle over a K3 surface;
\item[(c)] $q(\alpha)=-2$, $\alpha$ is integral, $Z$ is birationally a $\P^1$-bundle over the Hilbert square of a K3 surface;
\item[(d)] $q(\alpha)=-1/4$, $4\alpha$ is integral, $Z$ is birationally a $\P^1$-bundle over a product of two K3 surfaces;
\item[(e)] $q(\alpha)=-1$, $2\alpha$ is integral, $Z$ is birationally a $\P^1$-bundle over a product of two K3 surfaces.
\end{description}

It could be interesting to see whether one can make this more precise with the help of flops, as in the fourfold case, though there is no hope that one could suppress the word ``birationally'' altogether. Indeed the diagonal divisor (case (d)) is a contraction locus but not a $\P^1$-bundle, the fibers over the small diagonal are more complicated (cones over twisted cubics, \cite{T}).

\hfill

In the next section we perform a monodromy computation and use it to prove theorems \ref{_K3^3_MBM_Theorem_defo_} and 
\ref{_K3^2_MBM_Theorem_defo_}. In dimension 8 and 10, one derives a similar description of MBM classes and their loci.  \ref{_K3^3_MBM_Theorem_defo_} does not hold beyond dimension ten. Nevertheless one
can pursue the explicit description of extremal curves along these lines.

\section{Monodromy computations}
\label{_orbits_mono_Section_}



\theorem\label{_neg_def_mono_Theorem_}
Let $S$ be a K3 surface and let $\Lambda$ be the lattice $H^2(S^{[k]}, \Z)=H^2(S, \Z)\oplus \Z e$ with its BBF form.
If $k=2$, then for any primitive $z\in \Lambda$ with $q(z, z)<0$, there exists
$\gamma\in \tilde{O}(\Lambda)$ such that the rank 2 lattice
$\langle \gamma(z), e\rangle$ is negative definite,
unless $z$ and $e$ belong to the same  $\tilde{O}(\Lambda)$-orbit. If $3\leq k\leq 5$, there exists $\gamma\in \Mon$ such 
that the rank 2 lattice
$\langle \gamma(z), e\rangle$ is negative semidefinite.

\hfill

\remark In this statement $\tilde{O}(\Lambda)$ can be replaced by ``monodromy'', we just want to stress that the statement is about lattices. 

\hfill

This implies \ref{_K3^3_MBM_Theorem_defo_} and 
\ref{_K3^2_MBM_Theorem_defo_} because of the following proposition.

\hfill

\proposition Let $Z\subset X$ be a birational contraction
center on a deformation $X$ of the $n$-th Hilbert scheme of a K3 surface.
Suppose that \ref{_neg_def_mono_Theorem_} holds for its second cohomology lattice.
Then the pair $(X, Z)$ can be deformed to 
$(S^{[n]},Z')$, where $S$ is a non-algebraic K3 surface with (necessarily BBF nonpositive) cyclic Picard group.

\hfill

\proof Let $\eta$ be the MBM class associated with $Z$
and $\Teich_\eta$ the Teichm\"uller space of deformations
of $X$ such that $\eta$ remains of type (1,1). From the
global Torelli theorem it follows that
any lattice $\Lambda_1\subset \Lambda$ with
$\Lambda_1^\bot$ containing a positive 2-plane
and $\Lambda_1 \ni \eta$ can be realized as
a Picard lattice of $I\in \Teich_\eta$. Moreover all structures unseparables from $I$ (the points of $\Per^{-1}(\Per(I))$ have this Picard lattice, and the corresponding manifolds are all birational.
Applying this to the lattice $\langle \eta,\gamma^{-1}e\rangle$
from \ref{_neg_def_mono_Theorem_} , we obtain  deformations
of $(X, Z)$ with Picard lattice 
$\langle \eta,\gamma^{-1}e\rangle$. On such  deformations, $\eta$ remains of type $(1,1)$. Finally, some of them are Hilbert schemes 
of K3 surfaces: indeed the period points of Hilbert schemes of K3 surfaces in $\Perspace$ form the union of the orthogonal hyperplanes to the monodromy 
transforms of $e$ (each point of this union is a period point by Torelli theorem for K3 surfaces).
\endproof

\hfill

We denote by $\delta(z)$ the image of $\frac{z}{d(z)}$ 
in the discriminant group, where $d(z)$ denotes the divisibility of
$z$ (see \ref{orb-classif}). Recall from \ref{orb-classif}
that $z, z'\in \Lambda =(H^2(S^{[n]},\Z), q)$ belong to the same orbit of $\tilde{O}(\Lambda)$
if and only if $q(z) = q(z')$ and $\delta(z)=\delta(z')$.

\hfill

\example\label{gcd}
Let $S$ be a K3 surface, and $\Lambda =(H^2(S^{[2]},\Z), q)$,
so that $\Lambda= H^2(S,\Z)\oplus \Z e$,
where $q(e)=-2$, and the discriminant group of $\Lambda$ is  $\Z/2\Z$.
The primitive vectors $z\in \Lambda$ are $ax+be$, where
$x\in H^2(S,\Z)$ is primitive and $a, b\in \Z$ relatively prime. 
Since the lattice $H^2(S, \Z)$ is unimodular, $d(x)=1$ for any primitive
$x\in H^2(S, \Z)$. Therefore, $d(z)=2$ if and only if $a$ is even,
and $\delta(z)\neq 0$ if and only if $a$ is even and $b$ odd.
For such $z$, one has $q(z)=a^2q(x)- 2b^2$. More generally, if $\Lambda =(H^2(S^{[n]},\Z), q)$ and $e$ denotes half the diagonal class, 
then the divisibility of $ax+be$, $x\in H^2(S,\Z)$ primitive, is equal to g.c.d. $(a, 2(n-1)b)$.

\hfill

\noindent {\bf Proof of \ref{_neg_def_mono_Theorem_} in the case $k=2$}:
We consider two subcases,
$\delta(z)=0$ and $\delta(z)\neq 0$. If $\delta(z)=0$,
we can find $z'\in H^2(S,\Z)\subset H^2(S^{[2]},\Z)$
with $q(z)=q(z')$, and \ref{orb-classif}
implies that $z'=\gamma(z)$ for some $\gamma\in \tilde{O}(\Lambda)$.
Clearly, the lattice $\langle z', e\rangle$ is negative definite. 

Now, let $\delta(z)\neq 0$. Then $q(z)=a^2q(x)- 2b^2$, where $a=2 a_1$ is even, and
$q(z,z) = 4a_1^2 q(x) -2 b^2$ is equal to -2 modulo 8. Moreover,
any number $r= -2 \mod 8$ can be realized as $q(z')= 4q(x)-2$ for
$z'=2x'-e$, because any even number can be realized as a square of $x'\in H^2(S, \Z)$.
Since $q(z)< 0$, we have $q(z) \leq -10$ unless $q(z)=-2$.
In the latter case $z$ belongs to the same $\tilde{O}(\Lambda)$-orbit as $e$.
In the first case $z$ belongs to the same $\tilde{O}(\Lambda)$-orbit as
$z'=2x'-e$, with $q(x') \leq -2$. In this case,
the determinant of the lattice $\langle z', e\rangle$
is equal to $q(z') q(e) - q(z', e)^2 = -2 q(z') -4\geq 16$,
hence the lattice $\langle z', e\rangle$ is also negative definite. 
\endproof

\hfill

\noindent{\bf Proof of \ref{_neg_def_mono_Theorem_} in the case $k=3$}:
Now $e^2=-4$ and the discriminant group of the lattice is the cyclic group of order 4 
which we view as a subgroup of $\Q/\Z$,  and we have three possible subcases.

The first one $\delta(z)=0$ is the same as above.

If $\delta(z)=\pm{1/4}$, $z$ is conjugate to $e$ when $q(z)=-4$.
Otherwise we write $z$ as $4x+be$ where $x\in H^2(S,\Z)$ and $b$ is odd. 
Since $q(z)=16q(x)-4b^2$ is equal to -4 modulo 16 and 
is negative, we have $q(z) \leq -20$.
Moreover, $z$ is conjugate to some element of the form
$z'=4x'\pm e$, as $b^2-1$ is divisible by 8,
and hence $16q(x)-4b^2+4$ is divisible by 32. 
Now the lattice generated by $e$ and $z'$ is negative definite.

Some novelty appears when $\delta(z)=1/2$. Then one still finds a negative definite lattice whenever $q(z)\neq -4$,
but when $q(z)= -4$ then $z$ is conjugate to 
$z'$ of the form $2x'+e$ and the determinant of the lattice  
$\langle z', e\rangle$ is zero. \endproof

\hfill

\noindent{\bf Proof of \ref{_neg_def_mono_Theorem_} in the case of greater $k$}:
We set $t=2(k-1)$ so that $q(e)=-t$, and we eventually replace $z$ by a multiple in order to write uniformly $z=tx+be$, $b\in \Z$, $x$ primitive.
We first remark is that one may suppose that $|b|\leq t/2$. Indeed $tx+be$ has the same square as $tx'+(b\pm t)e$ for some 
$x'\in H^2(S,\Z)$ (the equation $t^2q(x)-tb^2=t^2q(x')-t(b^2\pm 2tb+t^2)$ on $q(x')$ has a solution since the terms $tb^2$ cancel and the other ones are divisible by
$2t^2$), and when $x'$ is primitive the two have the same
divisibility: indeed by \ref{gcd} this divisibility is g.c.d.$(t,bt)$ resp. g.c.d. $(t, (b\pm t)t)$, and the two are clearly equal. Therefore $tx+be$ and $tx'+(b\pm t)e$ are in the same monodromy orbit.

We know that $q(z)$ is negative. For small $t$ it means that $q(z)\leq -tb^2$: indeed this is the case when $q(x)$ is non-positive, and when $t\leq 8$ (that is, $k\leq 5$), already for $q(x)=2$ one has $2t^2-b^2t\geq 0$ provided $|b|\leq t/2$. Therefore the matrix of $q$ on $\langle z,e\rangle$ is negative semidefinite, and one concludes in the same way as before.

 Starting from $k=6$, the statement fails. Indeed, let $k=6$, so $t=10$, and $z=2x-e$ where $x$ is primitive of square 2, then $z^2=-2$ but the determinant of $q$ on $\langle z,e\rangle$ is negative and will remail negative when we 
replace $z$ by a conjugate. \endproof

\hfill

\remark In the last case mentioned, that is, $k=6$, the computation above shows that the only case when the pair $(M,z)$ cannot be deformed to the Hilbert scheme of a nonprojective surface with cyclic Picard group
is exactly as described: $z^2=-2$, $z=2x-e$ with $x$ orthogonal to $e$, $x$ of square $2$. On a generic deformation where
$e$ becomes a Hodge class, $M$ becomes the Hilbert scheme of a 
$K3$ surface $S$ with Picard group generated by $x$ of square $2$. Such an $S$ is a double covering of $\P^2$ and 
the MBM curve has class $z/2$. Up to monodromy, by \ref{fromcoverings} and the unicity statement (\ref{orb-classif}), it comes from a $g_4^1$ on a curve of genus 2 (the inverse image of a line under the covering map is of genus 2). Its MBM locus is birational to a $\P^2$-bundle over a product of a K3 surface with the symmetric square of a K3 surface.

This evokes the paper \cite{KLM} where many rational curves on  Hilbert schemes of a K3 surface and on generalized Kummer varieties are constructed from special linear systems on (possibly nodal) curves via Brill-Noether theory. Notice that in this particular case no Brill-Noether theory is needed, since obviously any curve of genus 2 carries a $g_4^1$. It seems likely and could be interesting to verify that such a description can be given in general.

\hfill

{\bf Acknowledgements:} We are grateful to Brendan Hassett, Ljudmila Kamenova, Emanuele Macri and Yuri Tschinkel for useful discussions.

\hfill

{

{\small
\noindent {\sc Ekaterina Amerik\\
{\sc Laboratory of Algebraic Geometry,\\
National Research University HSE,\\
Department of Mathematics, 7 Vavilova Str. Moscow, Russia,}\\
\tt  Ekaterina.Amerik@gmail.com}, also: \\
{\sc Universit\'e Paris-11,\\
Laboratoire de Math\'ematiques,\\
Campus d'Orsay, B\^atiment 425, 91405 Orsay, France}

\hfill

\noindent {\sc Misha Verbitsky\\
{\sc Instituto Nacional de Matem\'atica Pura e
              Aplicada (IMPA) \\ Estrada Dona Castorina, 110\\
Jardim Bot\^anico, CEP 22460-320\\
Rio de Janeiro, RJ - Brasil }\\
also:\\
{\sc Laboratory of Algebraic Geometry,\\
National Research University HSE,\\
Department of Mathematics, 7 Vavilova Str. Moscow, Russia,}\\
\tt  verbit@mccme.ru}.
 }
}

\end{document}